\magnification=\magstep1
\advance\voffset by 1 true cm
\baselineskip=15pt
\overfullrule =0pt 
\font\bbb=msbm10
\def\R{\hbox{\bbb R}}   
\def\N{\hbox{\bbb N}}    
\def\C{\hbox{\bbb C}}
\def\P{\hbox{\bbb P}}

\def\Z{\hbox{\bbb Z}}

\centerline {\bf HYPERSURFACE COMPLEMENTS, MILNOR FIBERS and}
\bigskip

 \centerline {\bf  MINIMALITY of ARRANGEMENTS  }

\vskip1.5truecm
\centerline {\bf by Alexandru Dimca }

\vskip1.5truecm

{\bf 1. The main results}

\bigskip

There is a gradient map associated to any reduced homogeneous polynomial $h \in \C[x_0,...,x_n]$ of degree $d>0$, namely
$$  grad(h): D(h) \to \P^n, ~~~~(x_0:...:x_n) \mapsto (h_0(x):...:h_n(x))$$
where $D(h)= \{x \in \P^n ; h(x) \not= 0\}$ is the principal open set associated to $h$ and $h_i={\partial h \over \partial x_i}$. Our first  result is the following topological description of the degree of the gradient map $  grad(h)$.

\bigskip

{\bf Theorem 1.} {\it For any reduced homogeneous polynomial   $h \in \C[x_0,...,x_n]$,
 the complement $D(h)$ is homotopy equivalent to a CW complex obtained from $D(h) \cap H$ by attaching $deg(grad(h)) $ cells of dimension $n$, where $H$ is a generic hyperplane
in $\P^n$. In particular, one has
$$deg(grad(h))= (-1)^n \chi (D(h) \setminus H).$$}
\bigskip
Note that the meaning of 'generic' here is quite explicit: the hyperplane $H$ has to be transversal to a stratification of the projective hypersurface $V(h)$ defined by $h=0$ in $\P^n$.

\bigskip

{\bf Corollary 2.} {\it Let $h^{(i)}$ denote the homogeneous polynomial obtained by restricting $h$ to a generic $i$-codimensional linear subspace in $\C^{n+1}$. Then
$$\chi (D(h))= \sum_{i=0,n} (-1)^{n-i}deg(grad(h^{(i)}))$$
where $deg(grad(h^{(n)})=1$ by convention.}

\bigskip

Using this result and the additivity of the Euler characteristic with respect to constructible partitions, one obtains formulas for the Euler characteristic
of any constructible set in terms of an alternating sum of degrees. This result should be compared with results by Szafraniec [Sz], where degrees of {\it real} polynomials play a similar role.

\bigskip

Let $f \in \C[x_0,...,x_n]$ be a homogeneous polynomial of degree $e>0$ with
global Milnor fiber $F=\{ x \in \C^{n+1}| f(x))=1\}$, see for instance [D1] for more on such varieties. Let $g: F \setminus N \to \R$
be the function $g(x)=h(x) {\overline h}(x)$, where $N=\{ x \in \C^{n+1}| h(x))=0\}$. Then we have the following.
\bigskip
{\bf Theorem 3.} {\it  For any reduced homogeneous polynomial   $h \in \C[x_0,...,x_n]$ and for any generic polynomial $f$ in the space of homogeneous polynomials of degree $e>0$  one has the following.

(i) the function $g$ is a Morse function.

(ii) the Milnor fiber $F$ is homotopy equivalent to a CW complex obtained from $F \cap N$ by attaching $|C(g)|$ cells of dimension $n$, where $C(g)$ is the critical set of the Morse function $g$.}

\bigskip

We point out that both Theorem 1. and Theorem 3. follow from the results by Hamm in [H]. In the case of Theorem 1. the homotopy type claim is a direct consequence from [H], Theorem 5. and also from Goresky and MacPherson [GM], Theorem 4.1, the new part being the relation between the number of $n$-cells and the degree of the gradient map $grad(h)$. We establish this equality by using polar curves, see section 2.

On the other hand, in Theorem 3. the main claim is that concerning the homotopy type and this follows from a very general result, see [H], Proposition 3. by a geometric argument described in section 3.

\bigskip

Our results above have interesting implications for the topology of hyperplane arrangements and these implications were our initial motivation in this study.
Let ${\cal A}$ be a hyperplane arrangement in the complex projective space $\P^n$, with $n>0$. Let $d>0$ be the number of hyperplanes in this arrangement and choose a linear equation $H_i: \ell_i(x)=0$ for each hyperplane $H_i$ in ${\cal A}$,
for $i=1,...,d.$

Consider the homogeneous polynomial $Q(x)= \prod _{i=1,d}\ell_i(x) \in \C[x_0,...,x_n]$ and the corresponding principal open set $M=D(Q)=\P^n \setminus \cup _{i=1,d}H_i$. The topology of the hyperplane arrangement complement $M$ is a central object of study in the theory of hyperplane arrangements, see Orlik-Terao [OT1].
As a  consequence of Theorem 1. we prove the following.

\bigskip

{\bf Corollary 4.}  {\it For any projective arrangement ${\cal A}$  as above one has
$$b_n(D(Q))=deg(grad(Q)).$$
In particular, the following are equivalent.

(i) the morphism $grad(Q)$ is dominant;

(ii)  $b_n(D(Q)) >0$;

(iii) the projective arrangement ${\cal A}$ is essential, i.e. the intersection $\cap _{i=1,d}H_i$ is empty.}

\bigskip

To obtain Corollary 4. from Theorem 1. all we need is the following.

\bigskip

{\bf Lemma 5.} {\it  For any arrangement ${\cal A}$  as above one has $(-1)^n
\chi (D(f) \setminus H)=b_n(D(f))$.}
\bigskip

This easy lemma has another very interesting consequence.
We say that a topological space $Z$ is minimal if $Z$ has the homotopy type of a CW-complex $K$ whose number of $k$-cells equals $b_k(K)$ for all $k \in \N$.

The importance of this notion for the topology of hyperplane arrangements was recently discovered by S. Papadima and A. Suciu, see [PS]  for various applications. The following result was independantly  obtained  by Randell, see [R], using similar techniques.

\bigskip

{\bf Corollary 6.} {\it The complement $M$ is a minimal space.}

\bigskip

It is easy to see that for $n>1$, the open set $D(f)$ is not minimal for $f$ generic of degree $d>1$ (just use $H_1(D(f),\Z)=\Z/d\Z$), but the Milnor fiber $F$ defined by $f$ is clearly minimal. Note that conversely, in spite of Corollary 6., the Milnor fiber $\{ Q=1\}$ associated to an arrangement is not  minimal in general.
\bigskip

From Theorem 3. we get a substantial strengthening of some of  the main results by Orlik and Terao in [OT2]. Let ${\cal A}'$ be the affine hyperplane arrangement in $\C^{n+1}$ associated to the projective arrangement ${\cal A}$. Note that $Q(x)=0$ is a reduced equation for the union $N$ of all the hyperplanes in ${\cal A}'$.
Let $f \in \C[x_0,...,x_n]$ be a homogeneous polynomial of degree $e>0$ with
global Milnor fiber $F=\{ x \in \C^{n+1}| f(x))=1\}$ and let $g: F \setminus N \to \R$
be the function $g(x)=Q(x) {\overline Q}(x)$ associated to the arrangement. The polynomial $f$ is called 
${\cal A}'$-generic if

(GEN1) the restriction of $f$ to any intersection $L$ of hyperplanes in ${\cal A}'$ is non-degenerate, in the sense that the associated projective hypersurface in $\P(L)$ is smooth, and

(GEN2) the function $g$ is a Morse function.

Orlik and Terao have shown in [OT2] that for an {\it essential} arrangement ${\cal A}'$, the set of ${\cal A}'$-generic functions $f$ is dense in the set of homogeneous polynomials of degree $e$, and, as soon as we have an ${\cal A}'$-generic function $f$, the following basic properties hold for any arrangement.

(P1) $~~b_q(F,F \cap N)=0$ for $q \not= n$ and

(P2) $~~b_{n}(F,F \cap N) \leq |C(g)|$, where $C(g)$ is the critical set of the Morse function $g$.

Moreover, for a special class of arrangements called {\it pure} arrangements it is shown in [OT2] that (P2) is actually an equality. In fact, the proof of (P2) in [OT2] uses some Morse theory, but we are unable to see the details behind the Corollary (3.5).

With this notation  the following is a direct consequence of Theorem 3.

\bigskip

{\bf Corollary 7.}

{\it For  any arrangement ${\cal A}'$ the following hold.

(i) the set of ${\cal A}'$-generic functions $f$ is dense in the set of homogeneous polynomials of degree $e>0$;

(ii)  the Milnor fiber $F$ is homotopy equivalent to a CW complex obtain from $F \cap N$ by attaching $|C(g)|$ cells of dimension $n$, where $C(g)$ is the critical set of the Morse function $g$. In particular $b_{n}(F,F \cap N)=|C(g)|$.  }

\bigskip

This paper represents a strengthening of the results in [D2] (in which the homological version of Theorem 1. and 3. above was proven).

The author thanks Stefan Papadima for raising the question answered by Corollary 4 above
and for lots of helpful comments.  In particular he informed me that Corollary 4 was proved by Paltin Ionescu in the case $n=2$ by completely different methods. I also thank Pierrette Cassou-Nogu\`es for drawing my attention on Richard Randell's preprint [R].

\vskip1.5truecm

{\bf 2. Polar curves, affine Lefschetz theory and degree of gradient maps}

\bigskip

The use of the local polar varieties in the study of singular spaces
is already a classical subject, see L\^e [L\^e],  L\^e -Teissier [LT] and the
references therein. 
Global polar curves in the study of the topology of polynomials (or, equivalently, the affine Lefschetz theory, for more on this equivalence see the beginning of the proof of Theorem 3.) is a topic under intense investigations, see for instance Cassou-Nogu\`es and Dimca [CD], Hamm [H], N\'emethi [N1-2], Siersma and Tib\u ar [ST], [T]. For all the proofs in this paper, the classical (local) theory is sufficient: indeed, all the objects being homogeneous, one can localize at the origin of $\C^{n+1}$ in the standard way, see [D1]. However, for the sake of geometric intuition, it seems to us easier to work with global (algebraic) objects, and hence we adopt this view-point in the sequel.

We recall briefly the notation and the results from [CD] and [N1-2]. Let $h \in \C[x_0, ... ,x_n]$ be a polynomial (even non-homogeneous to start with) and assume that the
fiber $F_t=h^{-1}(t)$ is smooth and connected, for some fixed $t \in \C$.

 For any hyperplane in $\P^n$, $ H:\ell =0$ where $  \ell(x)=h_0x_0+h_1x_1+...+h_nx_n $
we define the corresponding polar variety $\Gamma _H $ to be the union
of the irreducible components of the variety
$$ \{x \in \C^{n+1} \ |  \  \  rank(dh(x),d\ell (x))=1 \} $$
which are not contained in the critical set $S(h) = \{x \in \C^{n+1} \ | \
dh(x)=0 \} $ of $h$.

\bigskip

\noindent {\bf Lemma 8.} (see [CD], [ST])

{\it For a generic hyperplane $H$ we have the following properties.

(i) The polar variety $\Gamma _H$ is either empty or a curve,
i.e. each
irreducible component of $\Gamma _H$ has dimension 1.

(ii) dim$(F_t \cap \Gamma _H) \leq 0$ and the intersection
multiplicity
$(F_t,\Gamma_H)$ is independent of $H$.

(iii) The multiplicity $(F_t, \Gamma _H)$ is equal to the number of
tangent hyperplanes to $F_t$ parallel to the hyperplane $H$. For each
such tangent hyperplane $H_a$, the intersection $ F_t \cap H_a$ has
precisely one singularity, which is an ordinary double point.}

\bigskip

 The non-negative integer $(F_t, \Gamma _H)$ is called the polar
invariant
of the hypersurface $F_t$ 
and is denoted by $P(F_t)$. Note that $P(F_t)$ corresponds exactly to the classical notion of
class of a projective hypersurface, see [L].

We think of a projective hyperplane $H$ as above as the direction of an affine hyperplane $H' =\{ x \in \C^{n+1}|  \ell (x)=s \}$ for $s \in \C$. All the hyperplanes with the same direction form a pencil, and it is precisely the pencils of this type that are used in the affine Lefschetz theory, see [N1-2].
One of the  main results  in [CD] is the following, see also [ST] or [T] for similar results.
\bigskip

\noindent{\bf Proposition 9.}

\bigskip

{\it For a generic hyperplane $H'$ in the pencil of all hyperplanes in $\C^{n+1}$ with a fixed generic direction $H$, the fiber $F_t$ is homotopy equivalent to a CW-complex
obtained from  the section $F_t \cap H'$ by
attaching $P(F_t)$  cells of dimension $n$.
In particular}
$$P(F_t)=(-1)^{n}(\chi(F_t)- \chi(F_t \cap H'))=(-1)^n \chi(F_t \setminus H').$$
Moreover in this statement 'generic' means that the hyperplane $H'$ has to verify the following two conditions.

(g1) its direction, which is the hyperplane in $\P^n$ given by the homogeneous part of degree one in an equation for $H'$ has to be generic, and

(g2) the intersection $ F_t \cap H'$ has to be smooth.

These two conditions are not stated in [CD], but the reader should have no problem in checking them by using Theorem 3' in [CD] and the fact proved by N\'emethi in [N1-2] that the only bad sections in a good (i.e. the analog of a Lefchetz pencil in the projective Lefschetz theory, see [L]) pencil are the singular sections. Completely similar results hold for generic pencils with respect to a closed smooth subvariety $Y$ in some affine space $\C^N$, see [N1-2], but note that the polar curves are not mentionned there.
\bigskip

{\bf Proof of Theorem 1.}
\bigskip

In view of Hamm's affine  Lefschetz theory, see [H], Theorem 5. and also from Goresky and MacPherson [GM], Theorem 4.1 , the only thing to prove is the equality between the number $k_n$ of $n$-cells attached and the degree of the gradient.

Assume from now on that the polynomial $h$ is homogeneous of degree $d$ and that $t=1$. It follows from (g1) and (g2) above that we may choose the generic hyperplane $H'$ passing through the origin.

Moreover, in this case, the polar curve $\Gamma _H$, being defined by homogeneous equations, is a union of lines $L_j$ passing through the origin. For each such line we choose a parametrization $t \mapsto a_j t$ for some $a_j \in \C^{n+1}, a_j \not=0$. It is easy to see that the intersection $F_1 \cap L_j$ is either empty (if $h(a_j)=0$) or consists of exactly $d$ distinct points with multiplicity one (if $h(a_j)\not=0$). The lines of the second type are in bijection with the points in $grad(h)^{-1}(D_{H'})$,
where $D_{H'}\in \P^n$ is the point corresponding to the direction of the hyperplane $H'$. It follows that
$$ d \cdot deg(grad(h))= P(F_1).$$

The $d$-sheeted unramified coverings $ F_1 \to D(h)$ and $F_1 \cap H' \to D(h) \cap H$ give the result, where $H$ is the projective hyperplane corresponding to the affine hyperplane (passing through the origin) $H'$. Indeed, they imply the equalities: $\chi(F_1)=d \cdot \chi (D(h))$ and $\chi(F_1 \cap H')=d \cdot \chi (D(h) \cap H)$. Hence we have $deg(grad(h))= \chi (F_1,F_1 \cap H')/d=\chi(D(h),D(h) \cap H)=k_n.$

\bigskip

{\bf Remark 10.} The gradient map $grad(h)$ has a natural extension to the larger open set $D'(h)$ where at least one of the partial derivatives of $h$ does not vanish. It is obvious (by a dimension argument) that this extension has the same degree as the map $grad(h)$.

\vskip1.5truecm

{\bf 3. Non-proper Morse Theory}

\bigskip

For the convenience of the reader we recall, in the special case we need, a basic result of Hamm, see [H], Proposition 3, with our addition  concerning the condition $(c0)$ in [DP], see Lemma 3. and Example 2. The final claim on the number of cells to be attached is also standard, see for instance [ST] and [T].

\bigskip

\noindent {\bf Proposition 11.}

{\it Let $A$ be a smooth algebraic subvariety in $\C^p$ with dim$A=m$. Let $f_1,...,f_p $ be polynomials in $\C[x_1,...,x_p]$.
For $ 1 \leq j \leq p$, denote by  $\Sigma _j$ the set of critical points of the mapping $(f_1,...,f_j): A \setminus \{z \in A; f_1(z)=0 \} \to \C^j$ and let $\Sigma _j'$ denote the closure of $\Sigma _j$ in $A$.
Assume  that the following conditions hold.

(c0) The set $\{z \in A; |f_1(z)| \leq a_1, ..., |f_p(z)| \leq a_p \}$ is compact for any positive numbers $a_j$, $j=1,...,p$.

(c1) The  critical set $\Sigma _1$ is finite.

(cj) (for $j=2,...,p$) The map  $(f_1,...,f_{j-1}): \Sigma _j' \to \C ^{j-1}$ is  proper.

Then $A$ has the homotopy type of a space obtained from $A_1=\{x \in A; f_1(x)=0\}$ by attaching $m$-cells and the number of these cells is the sum of the Milnor numbers $\mu(f_1,x)$ for $x \in \Sigma _1$.}

\bigskip

{\bf Proof of Theorem 3.}

\bigskip

We set $X=h^{-1}(1)$.
Let  $v: \C^{n+1} \to \C^N$ be the Veronese mapping of degree $e$ sending $x$ to all the monomials of degree $e$ in $x$ and set $Y=v(X)$.
Then $Y$ is a smooth closed subvariety in $\C^N$ and $v:X \to Y$ is an unramified (even Galois) covering of degree $c$, where $c=g.c.d.(d,e)$. To see this, use the fact that $v$ is a closed immersion on $\C^N \setminus \{0\}$ and $v(x)=v(x')$ iff $x'=u \cdot x$ with $u^c=1$.

Let $H$ be a generic hyperplane direction in $\C^N$ with respect to the subvariety $Y$ and let $C(H)$ be  the finite set of all the points $p \in Y$ such that there is an affine hyperplane $H'_p$ in the pencil determined by $H$ that is tangent to $Y$ at the point $p$ and the intersection $Y \cap H'_p$ has a complex Morse (alias non-degenerated, alias $A_1$) singularity. 
Under the Veronese mapping $v$, the generic hyperplane direction $H$ corresponds to a homogeneous polynomial of degree $e$ which we call from now on $f$.

To prove the first claim (i) we proceed as follows. It is known that doing affine Lefschetz theory for a pencil of hypersurfaces $\{h=t\}$ is equivalent to doing (non-proper) Morse theory for the function $|h|$ or, what amounts to the same, for the function $|h|^2$. More explicitly, in view of the last statement at the end of the proof of Lemma (2.5) in [OT2] (which clearly applies to our more general setting since all the computations there are local), $g$ is a Morse function iff each critical point of $h:F \setminus N \to \C$ is an $A_1$-singularity. Using the homogeneity of both $f$ and $h$, this last condition on $h$ is equivalent to the fact that each critical point of the function $f:X \to \C$ is an $A_1$ singularity, condition fulfilled in view of the choice of $H$ and since $v:X \to Y$ is a local isomorphism.

Now we pass on to the proof of the claim (ii) in Theorem 3. Recall first that any polynomial function $h:\C^{n+1} \to \C$ admits a Whitney stratification satisfying Thom $a_h$-condition: this is a constructible stratification ${\cal S}$ such that the open stratum, say $S_0$, coincides with the set of regular points for $h$ and for any other stratum, say $S_1 \subset h^{-1}(0)$, and any sequence of points $q_m \in S_0$ converging to $q \in S_1$ such that the sequence of tangent spaces $T_{q_m}(h)$ has a limit $T$, then $T_qS_1 \subset T$, see Hironaka [Hi], Corollary 1, page 248 (and note that the requirement of $f$ proper in that Corollary is not necessary in our case, as any algebraic map can be compactified). Here and in the sequel, for a map $\phi : {\cal X} \to {\cal Y}$ and a point $q \in {\cal X}$ we denote by $T_q(\phi)$ the tangent space to the fiber $\phi ^{-1}(\phi (q))$ at the point $q$, assumed to be a smooth point on this fiber.

Since in our case $h$ is a homogeneous polynomial, we can find a stratification ${\cal S}$ as above such that all of its strata are $\C^*$-invariant, with respect to the natural $\C^*$-action on $\C^{n+1}$. In this way we obtain an induced Whitney stratification  ${\cal S'}$ on the projective hypersurface $V(h)$.
We choose our polynomial $f$ such that the corresponding projective hypersurface $V(f)$ is smooth and transversal to the stratification  ${\cal S'}$. In this way we get  an induced Whitney stratification  ${\cal S'}_1$ on the projective complete intersection $V_1= V(h) \cap V(f)$.

We use Proposition 11. above with $A=F$ and $f_1=h$. All we have to show is the existence of polynomials $f_2,...,f_{n+1}$ satisfying the conditions listed in Proposition 11.

We will choose these polynomials inductively to be generic linear forms as follows. We choose $f_2$ such that the corresponding hyperplane $H_2$ is transversal to the stratification  ${\cal S'}_1$. Let ${\cal S'}_2$ denote the induced stratification on $V_2 =V_1 \cap H_2$. Assume that we have constructed $f_2,...,f_{j-1}$, ${\cal S'}_1,...,{\cal S'}_{j-1}$ and $V_1,...,V_{j-1}$. We choose $f_j$ such that the corresponding hyperplane $H_j$ is transversal to the stratification  ${\cal S'}_{j-1}$. Let ${\cal S'}_j$ denote the induced stratification on $V_j =V_{j-1} \cap H_j$. Do this for $j=3,...,n$ and choose for $f_{n+1}$ any linear form.

With this choice it is clear that for $1 \leq j \leq n$, $V_j$ is a complete intersection of dimension $n-1-j$. In particular, $V_n=\emptyset$, i.e.

$(c0')~~~~~~~\{x \in \C^{n+1}; f(x)=h(x)=f_2(x)=...=f_{n}(x)=0\}=\{0\}.$

Then the map $(f,h,f_2,...,f_n):\C^{n+1} \to \C^{n+1}$ is proper, which clearly implies the condition $(c0)$.

The condition $(c1)$ is fulfilled by our construction of $f$. Assume that we have already checked that the conditions $(ck)$ are fulfilled for $k=1,...,j-1$. We explain now why  the next condition $(cj)$ is fulfilled.

Assume that the condition $(cj)$ fails. This is equivalent to the existence of a sequence $p_m$ of points in $\Sigma _j'$ such that

$(*)~~~~~|p_m| \to \infty$ and $f_k(p_m) \to b_k$ (finite limits) for $1 \leq k \leq j-1$.

Since $\Sigma _j$ is dense in
$\Sigma _j'$, we can even assume that $p_m \in  \Sigma _j$.

Note that $\Sigma _{j-1} \subset \Sigma _j$ and the condition $c(j-1)$ is fulfilled. This implies that we may choose our sequence $p_m$ in the difference
$\Sigma _j \setminus  \Sigma _{j-1}$. In this case we get
$$(**)~~~~~ f_j \in Span(df(p_m),dh(p_m),f_2,...,f_{j-1})$$
the latter being a $j$-dimensional vector space.

Let $q_m ={p_m \over |p_m|} \in S^{2n+1}$. Since the sphere $S^{2n+1}$ is compact we can assume that the sequence $q_m$ converges to a limit point $q$. By passing to the limit in $(*)$ we get $q \in V_{j-1}$. Moreover, we can assume (by passing to a subsequence) that the sequence of $(n-j+1)$-planes $T_{q_m}(h,f,f_2,...,f_{j-1})$ has a limit $T$. Since $p_m \notin \Sigma_{j-1}$, we have
$$ T_{q_m}(h,f,f_2,...,f_{j-1})=T_{q_m}(h) \cap T_{q_m}(f) \cap H_2 \cap ... \cap H_{j-1}$$
As above, we can assume that the sequence $T_{q_m}(h)$ has a limit $T_1$ and, using the $a_h$-condition for the stratification ${\cal S}$ we get $T_qS_i \subset T_1$ if $q \in S_i$. Note that we have $T_{q_m}(f) \to T_{q}(f)$ and hence
$T=T_1 \cap T_{q}(f) \cap H_2 \cap ... \cap H_{j-1}$. It follows that
$$T_qS_{i,j-1}= T_qS_i \cap T_{q}(f) \cap H_2 \cap ... \cap H_{j-1} \subset T$$
where $S_{i,j-1}= S_i \cap V(f) \cap H_2 \cap ... \cap H_{j-1}$ is the stratum corresponding to the stratum $S_i$ in the stratification ${\cal S'}_{j-1}$. On the other hand, the condition $(**)$ implies that $T_qS_{i,j-1} \subset T \subset   H_j $, a contradiction to the fact that $H_j$ is transversal to ${\cal S'}_{j-1}$.

\vskip1.5truecm

{\bf 4. Complements of hyperplane arrangements}

\bigskip

{\bf Proof of Lemma 5}.

Here we just give the main idea, since the details are standard. One has to use the method of deletion and restriction, see [OT1], p. 17, the obvious additivity of the Euler characteristics  and, more subtly, the additivity of the top Betti numbers coming from the exact sequence (8) in [OT1], p. 20 or (3.8) in [DL].

\bigskip
{\bf Proof of Corollary 4}.

To complete this proof  we only have to explain why the claims $(ii)$ and $(iii)$ are equivalent. If the projective arrangement is not essential, then using a projection onto $\P^{n-1}$ with center a point in all the hyperplanes $H_i$ we get a fiber bundle $D(Q) \to U$ with fiber $\C$ and base $U$, an affine variety of dimension $n-1$. This implies $b_n(D(Q))=0$.

If the arrangement is essential, then $d \geq n+1$ and we may assume that
$\ell _i(x)=x_{i-1}$ for $i=1,...,n+1.$ In the case $d = n+1$, we are done, since in this case $D(Q)= (\C^*)^n$ and hence $b_n(D(Q))=1$. In the remaining case 
$d>n+1$, one should use the additivity of the top Betti numbers alluded above in the proof of Lemma 5.
\bigskip

{\bf  Proof of Corollary 6.}

\bigskip

Using the Affine Lefschetz Theorem of Hamm, see Theorem 5 in [H], we know that for a generic projective hyperplane $H$, the space $M$ has the homotopy type of a space obtained from $M \cap H$ by attaching $n$-cells. The number of these cells is given by
$$(-1)^n \chi (M, M \cap H)=(-1)^n \chi (M \setminus H)=b_n(M)$$
see Lemma 4. above.

To finish the proof of the minimality of $M$ we proceed by induction using the equalities
$$b_k(M)=b_k(M \cap H)$$
for $0 \leq k <n$. Indeed, for $0 \leq k <n-1$, this is obvious since we attach only $n$-cells. The equality for $k=n-1$ follows from these previous equalities and a new application of Lemma 5.

\bigskip

{\bf Remark 12.}

Let $\mu_e$ be the cyclic group of the $e$-roots of unity. Then there is a natural algebraic action of $\mu_e$ on the space $F \setminus N$ occuring in Theorem 1'. The corresponding weight equivariant Euler polynomial (see [DL] for a definition) gives information on the relation between the induced $\mu_e$-action on the cohomology $H^*(F \setminus N)$ and the functorial Deligne mixed Hodge structure present on cohomology.

When $N$ is a hyperplane arrangement ${\cal A}'$  and $f$ is an ${\cal A}'$-generic function , this weight equivariant Euler polynomial
can be combinatorically computed from the lattice associated to the arrangement
(see Corollary (2.3) and Remark (2.7) in [DL]) using the fact that the 
 weight equivariant Euler polynomial of the $\mu_e$-variety $F$ is known, see for instance [D1].

\bigskip

\vskip1.5truecm

\noindent {\bf REFERENCES}
\bigskip

\item{[CD]} Pi. Cassou-Nogu\`es, A. Dimca: Topology of complex polynomials via polar curves, Kodai Math. J. 22(1999), 131-139.

\item{[D1]} A. Dimca: Singularities and Topology of Hypersurfaces,
Universitext,Springer, 1992.

\item{[D2]} A. Dimca: Arrangements, Milnor fibers and polar curves, math.AG/0011073.

\item {[DL]} A. Dimca, G.I. Lehrer: Purity and equivariant weight
polynomials, dans le volume: Algebraic Groups and Lie Groups, editor
G.I. Lehrer, Cambridge University Press, 1997.

\item{[DP]} A. Dimca, L. P\u aunescu: On the connectivity of complex affine hypersurfaces,II, Topology 39 (2000),1035-1043.

\item{[GM]} M. Goresky, R. MacPherson: Stratified Morse theory,  Proc. Symp. Pure Math., Singularities, Volume 40, Part 1 (1983), 517-533.

\item{[H]} H. A. Hamm: Lefschetz theorems for singular varieties, Proc. Symp. Pure Math., Singularities, Volume 40, Part 1 (1983), 547-557.

\item{[Hi]} H. Hironaka: Stratifications and flatness,in: Real and Complex Singularities, Sijthoff and Noordhoff, 1977, 199-265.

\item{[L]} K. Lamotke: The topology of complex projective varieties
after S. Lefschetz, Topology 20(1981),15-51.

\item {[L\^e]} D.T. L\^e: Calcul du nombre de cycles \'evanouissants d'une hypersurface complexe, Ann. Inst. Fourier, Grenoble 23 (1973), 261-270.

\item{[LT]} D.T. L\^e, B. Teissier: Vari\'et\'es polaires locales et
classes de Chern des vari\'et\'es singuli\`eres, Ann. Math. 114 (1981),457-491.

\item{[N1]} A. N\'emethi: Th\'eorie de Lefschetz pour les vari\'et\'es
alg\'ebriques affines, C. R. Acad. Sci. Paris 303(1986), 567-570.

\item{[N2]} A. N\'emethi: Lefschetz theory for complex affine
varieties, Rev. Roum. Math. Pures et Appl. 33 (1988),233-260.

\item{[OT1]} P. Orlik, H. Terao: Arrangements of hyperplanes, Springer 1992.

\item{[OT2]} P. Orlik, H. Terao: Arrangements and Milnor fibers, Math. Ann. 301
(1995), 211-235.

\item{[PS]} S. Papadima, A. Suciu: Higher homotopy groups of complements of complex hyperplane arrangements, math.AT/0002251.

\item{[R]} R. Randell: Morse theory, Milnor fibers and hyperplane arrangements, math.AT 0011101.

\item{[ST]} D. Siersma, M. Tib\u ar: Singularities at infinity and their vanishing cycles II, Monodromy, Publ. RIMS, to appear.

\item{[Sz]} Z. Szafraniec: On the Euler characteristic of complex algebraic varieties, Math. Ann. 280 (1988),177-183.

\item{[T]} M. Tib\u ar: Asymptotic equisingularity and topology of complex hypersurfaces, Int. Math. Res. Not. 18 (1998), 979-990.

\bigskip

Laboratoire de Math\'ematiques Pures de Bordeaux

Universit\'e Bordeaux I

33405 Talence Cedex, FRANCE

\bigskip

email: dimca@math.u-bordeaux.fr

\bye